\documentclass[12pt,twoside,reqno]{amsart}
\usepackage[colorlinks=true,citecolor=blue]{hyperref}
\usepackage{mathptmx, amsmath, amssymb, amsfonts, amsthm, mathptmx, enumerate, color}
\usepackage{graphicx}
\usepackage{epstopdf}

\newtheorem{theorem}{Theorem}[section]
\newtheorem{corollary}{Corollary}[section]
\newtheorem{lemma}{Lemma}[section]

\theoremstyle{definition}

\newtheorem{conjecture}{Conjecture}

\numberwithin{equation}{section}

\begin{document}
\setcounter{page}{1}

\vspace*{1.0cm}
\title[Norm of the Hilbert matrix operator on logarithmically weighted Bloch  and Hardy spaces]
{Norm of the Hilbert matrix operator on logarithmically weighted Bloch  and Hardy spaces}
\author[ S. Ye, Q. Zheng]{ Shanli Ye$^{*}$, Qisong Zheng}
\maketitle
\vspace*{-0.6cm}

\begin{center}
{\footnotesize {\it

School of Science, Zhejiang University of Science and Technology, Hangzhou 310023, China.

}}\end{center}

\vskip 4mm {\small\noindent {\bf Abstract.}
In this paper, we compute the exact value of the norm of the Hilbert matrix operator $\mathcal{H}$ acting from the classical Bloch space  $\mathcal{B}$  into the logarithmically weighted Bloch space  $\mathcal{B}_{\log}$, and show that it equals $\frac{3}{2}$; we also find that the norm from the space of bounded analytic functions $H^\infty$  into  the logarithmically weighted Hardy space  $H^{\infty}_{\log}$ is  $1$. Furthermore, we establish both lower and upper bounds for the norm of $\mathcal{H}$ when it maps  from the $\alpha$-Bloch space $\mathcal{B}^\alpha$ into the logarithmically weighted $\mathcal{B}^\alpha_{\log}$ with $1 <\alpha < 2$, and from the Hardy space $H^{1}$ into the logarithmically weighted Hardy space $H^{1}_{\log}$.

 }

\renewcommand{\thefootnote}{}
\footnotetext{ $^*$Corresponding author.
\par
E-mail addresses:  slye@zust.edu.cn (S. Ye).
\par
}

\section{Introduction}
The Hilbert matrix operator $\mathcal{H}$ is a fundamental object in operator theory. In recent years, significant research activity has been focused on analyzing its boundedness, norm, and other properties on various spaces of analytic functions. The investigation into its boundedness and norm, in particular, has been a very active field of study \cite{Ale,Bra,Bou,Bre,Dai1,Dai2,Lan}.
The initial work by Diamantopoulos and Siskakis \cite{Dia2} established the boundedness of $\mathcal{H}$ on the Hardy space $H^p$ for $1<p<\infty$, concurrently providing an initial upper bound for its norm. This research was later extended by Diamantopoulos \cite{Dia1} to the Bergman spaces $A^p$ for $2<p<\infty$, where a corresponding upper bound for the norm of $\mathcal{H}$ was also determined. Building on this work, Dostani\'{c} et al. \cite{Dos} obtained the exact norm of $\mathcal{H}$ on the Hardy space $H^p$ ($1<p<\infty$) and provided its precise norm on the Bergman space $A^p$ for $4 \leq p < \infty$. In 2018, Bo\v{z}in and Karapetrovi\'{c} \cite{Bou} derived the exact norm of $\mathcal{H}$ on the Bergman space $A^p$ for $2 < p < 4$. The investigation into the boundedness of $\mathcal{H}$ on weighted Bergman spaces $A^p_\alpha$ began with \cite{Gal}. This initial work spurred extensive research into the norms of $\mathcal{H}$ on $A^p_\alpha$ across different $\alpha$, as detailed in \cite{BTW,Dai1,DK,Jev1,Jev2,Kar1,Kar2,Lan,Lin1,Lin2,N}. In 2023, the first author and Feng \cite{Ye} applied Littlewood's classical result \cite[pp. 93-96]{Lit} to estimate the norm of $\mathcal{H}$ from the logarithmically weighted space $A_{\log^\alpha}^2$ to $A^2$ for $\alpha>2$, extending on the work by Jevti\'{c} and Karapetrovi\'{c} \cite{Jev1}, and Karapetrovi\'{c} \cite{Kar3}. Very recently, Hu  and the first author\cite{Hu} computed the exact norm of the Hilbert matrix operator $\mathcal{H}$ acting from the logarithmically weighted Korenblum space $H^\infty_{\alpha,\log}$ into the Korenblum space $H^\infty_\alpha$, and from the Hardy space $H^\infty$ into the classical Bloch space $\mathcal{B}$. It was  shown in \cite{Lan} that $\mathcal{H}$ is unbounded on the  Bloch space $\mathcal{B}$.  Moreover, it is known from \cite{Dia2} that $\mathcal{H}$ is unbounded on the Hardy spaces $H^p$ for $p=1$ and $p=\infty$.

   In this paper, we introduce the logarithmically weighted Hardy spaces $H^p_{\log}$ and the logarithmically weighted Bloch spaces $\mathcal{B}^\alpha_{\log}$. We compute the exact values of the norms of the Hilbert matrix operator $\mathcal{H}$: specifically, it is shown to be $1$ when $\mathcal{H}$ acts from the space of bounded analytic functions $H^\infty$ into the logarithmically weighted Hardy space $H^{\infty}_{\log}$, and $\displaystyle\frac{3}{2}$ when it acts from the classical Bloch space $\mathcal{B}$ into the logarithmically weighted Bloch space $\mathcal{B}_{\log}$. Furthermore, we establish both lower and upper bounds for the norm of $\mathcal{H}$ when it maps from the Hardy space $H^{1}$ into the logarithmically weighted Hardy space $H^{1}_{\log}$, and from the $\alpha$-Bloch space $\mathcal{B}^\alpha$ into the logarithmically weighted Bloch space $\mathcal{B}^\alpha_{\log}$ for $1 < \alpha < 2$.

\section{Preliminaries}

Let $\mathbb{D}$ denote the open unit disk of the complex plane $\mathbb{C}$. Additionally, let $H(\mathbb{D})$ denote the set of all analytic functions in  $\mathbb{D}$.

For $0 < p \leq \infty $, the Hardy space $H^p$ is the space comprising all functions  $f \in H(\mathbb{D} )$ such that
$$\|f\|_{H^p}=\sup_{0\leq r <1} M_p(r,f)<\infty,$$
where
$$M_p(r,f)=\left( \frac{1}{2\pi}\int_0^{2\pi} |f(re^{it})|^p dt\right)^{\frac{1}{p}}, \quad 0<p<\infty;$$
$$M_\infty(r,f)=\sup_{0\leq t<2\pi}|f(re^{it})|.$$
References such as \cite{Dur1} are suggested for further details on the notation and results about Hardy spaces.

Now  we define the logarithmically weighted Hardy space $H^p_{\log}$, which consists of those $f \in H(\mathbb{D})$ such that
$$\|f\|_{H^p_{\log}}\overset{def}{=} \sup_{0\leq r <1}\left(\log\frac{e}{(1-r)^2}\right)^{-1} M_p(r,f)<\infty.$$
  It is easy to see that $H^p \subseteq H^p_{\log}$.

For $0<\alpha<\infty$, the $\alpha$-Bloch space $\mathcal{B}^\alpha$ consists of those functions $f \in H(\mathbb{D})$ with
$$\|f\|_{\mathcal{B}^\alpha}=|f(0)|+\|f\|_{\alpha*}<\infty,$$
where
$$\|f\|_{\alpha*}=\sup_{z\in\mathbb{D}}(1-|z|^2)^\alpha|f'(z)|.$$
Note that $\mathcal{B}^1$ is the classical Bloch space $\mathcal{B}$, and we write  $\|f\|_*=\|f\|_{1*}$ for brevity.  We refer to \cite{And, Zhu} as general references for the classical Bloch space and the $\alpha$-Bloch spaces.

We introduce, for $0<\alpha<\infty$, the logarithmically weighted $\alpha$-Bloch space, denoted by $\mathcal{B}^\alpha_{\log}$, as the set of all
 $f \in H(\mathbb{D})$ satisfying
$$\|f\|_{\mathcal{B}^\alpha_{\log}}\overset{def}{=}|f(0)|+\sup_{z\in\mathbb{D}}\left(\log\frac{e}{(1-|z|)^2}\right)^{-1}(1-|z|^2)^\alpha|f'(z)|<\infty.$$
We write $\mathcal{B}_{\log}$  for the case  $\alpha=1$, and it follows that  $\mathcal{B}^\alpha \subseteq \mathcal{B}^\alpha_{\log}.$

The Hilbert matrix is an infinite matrix $\mathcal{H}$ whose entries are $a_{n,k} = \frac{1}{n+k+1}$, $n,k \geq 0$. The Hilbert matrix $\mathcal{H}$ can also be viewed as an operator on spaces of analytic functions by its action on their Taylor coefficients. Hence for $f \in H(\mathbb{D})$ with Taylor expansion $f(z) = \sum_{k=0}^\infty a_k z^k$, we define a transformation $\mathcal{H}$ by
\begin{align}
	\mathcal{H}f(z)=\sum_{n=0}^\infty \left( \sum_{k=0}^\infty \frac{a_k}{n+k+1} \right)z^n,\notag
\end{align}
whenever the right hand side makes sense and defines an analytic function in $\mathbb{D}$. Consequently, we can assert that
\begin{align} \label{eq1.1}
	\mathcal{H}f(z) &= \sum_{n=0}^\infty \left( \sum_{k=0}^\infty \frac{a_k}{n+k+1} \right)z^n  \notag\\
	&= \sum_{n=0}^\infty \left( \sum_{k=0}^\infty a_k \int_0^1 t^{n+k}dt \right)z^n
	\notag\\
	&= \int_0^1 \sum_{k=0}^\infty a_k  t^k \sum_{n=0}^\infty  t^n z^n dt \notag \\
	&= \int_0^1 \frac{f(t)}{1-tz} dt.
\end{align}
By (\ref{eq1.1}), we change the path of integration to obtain the weighted composition operators $T_t$ (see \cite{Dia2}):
\begin{align*}
	\mathcal{H} f(z)=\int_0^1 T_t f(z) dt,
\end{align*}
where
\begin{align}
	T_tf(z)=w_t(z)f(\phi_t(z)),  \quad w_t(z)=\frac{1}{1-(1-t)z}, \quad \phi_t(z)=\frac{t}{1-(1-t)z}. \notag
\end{align}
According to (\ref{eq1.1}), we obtain
\begin{align}\label{eq3.1}
	(\mathcal{H}f)'(z)=\int_{0}^{1}\frac{t}{(1-tz)^2}f(t)dt.
\end{align}
For $z\in\mathbb{D}$, we can choose the path in \cite{Dia2}$$\zeta(t)=\zeta_z(t)=\frac{t}{(t-1)z+1}.$$
The change of variable in (\ref{eq3.1}) gives
\begin{align}\label{eq3.2}
	(\mathcal{H}f)'(z)=\int_{0}^{1}\frac{t}{\left[(t-1)z+1\right](1-z)}f(\phi_t(z))dt.
\end{align}
For $f\in\mathcal{B}^\alpha$ and $z\in\mathbb{D}$. We have
	$$
		|f(z)-f(tz)|=\left|z\int_{t}^{1}f'(sz)ds\right|\leq\|f\|_{\alpha*}\int_{t}^{1}\frac{|z|}{(1-|sz|^2)^{\alpha}}ds\leq\|f\|_{\alpha*}\int_{|z|t}^{|z|}\frac{1}{(1-x^2)^{\alpha}}dx.
	$$
	In particular, letting $t=0$, we have
\begin{align}\label{eq3.4}
	|f(z)|\leq\|f\|_{\alpha*}\int_{0}^{|z|}\frac{1}{(1-x^2)^{\alpha}}dx+|f(0)|
	\begin{cases}
		&=\displaystyle\frac{1}{2}\log\frac{1+|z|}{1-|z|}\|f\|_{\alpha*}+|f(0)|\quad if\quad\alpha =1,\\
		&\leq \displaystyle\frac{(1-|z|)^{1-\alpha}-1}{\alpha-1}\|f\|_{\alpha*}+|f(0)|\quad if\quad\alpha \not=1.\\
	\end{cases}
\end{align}

\section{The norm of the Hilbert matrix $\|\mathcal{H}\|_{\mathcal{B} \rightarrow \mathcal{B}_{\log}}$}\label{s5}
In \cite[p.171]{Lan} , it is shown  that $\mathcal{H}$ is not bounded on the classical Bloch space $\mathcal{B}$.  However, we establish that $\mathcal{H}$ is bounded from $\mathcal{B}$ to $\mathcal{B}_{\log}$, and compute the exact norm  $\|\mathcal{H}\|_{\mathcal{B}\rightarrow \mathcal{B}_{\log}}$ in this section.

\begin{lemma}\label{Le3.1}
	The norm of the Hilbert matrix operator acting from $\mathcal{B}$ into $\mathcal{B}_{\log}$ satisfies
	\begin{align*}
		\| \mathcal{H}\|_{\mathcal{B}\rightarrow\mathcal{B}_{\log}}= max\{1+\sup_{0\leq r<1}(1+r)\left(\log\frac{e}{(1-r)^2}\right)^{-1}\int_{0}^{1}\frac{t}{1+(t-1)r}dt, \\ \log2+\frac{1}{2}\sup_{0\leq r<1}(1+r)\left(\log\frac{e}{(1-r)^2}\right)^{-1}
		\int_{0}^{1}\frac{t}{1+(t-1)r}\log\frac{1+t+(t-1)r}{(1-t)(1-r)}dt\}.
	\end{align*}
	
\end{lemma}

\begin{proof}
		Suppose $f\in\mathcal{B}$ and $z\in\mathbb{D}$. 	By (\ref{eq3.2}) and (\ref{eq3.4}), we obtain that
	\begin{align*}
    \|\mathcal{H}f\|_{\mathcal{B}_{\log}}&=|\mathcal{H}f(0)|+\sup_{z\in\mathbb{D}}(1-|z|^2)\left(\log\frac{e}{(1-|z|)^2}\right)^{-1}|(\mathcal{H}f)'(z)|\\
        &\leq\int_{0}^{1}\left|f(t)\right|dt+\sup_{z\in\mathbb{D}}(1-|z|^2)\left(\log\frac{e}{(1-|z|)^2}\right)^{-1}\int_{0}^{1}\left|\frac{\phi_t(z)}{1-z}f(\phi_t(z))\right|dt\\
		&\leq\|f\|_{*}\int_{0}^{1}\frac{1}{2}\log\frac{1+t}{1-t}dt+|f(0)|+\sup_{z\in\mathbb{D}}(1+|z|)\left(\log\frac{e}{(1-|z|)^2}\right)^{-1}\\
        &\times\int_{0}^{1}\left(\frac{1}{2}\log\frac{1+|\phi_t(z)|}{1-|\phi_t(z)|}\left|\phi_t(z)\right|\|f\|_{*}+|\phi_t(z)||f(0)|\right)dt.\\
	\end{align*}
    Since $|\phi_t(z)| \leq \phi_t(|z|)$,  we obtain
    \begin{align*}	
    &\|\mathcal{H}f\|_{\mathcal{B}_{\log}}\leq\left(1+\sup_{0\leq r<1}(1+r)\left(\log\frac{e}{(1-r)^2}\right)^{-1}\int_{0}^{1}\frac{t}{1+(t-1)r}dt\right)|f(0)|\\
		&+\left(\log2+\frac{1}{2}\sup_{0\leq r<1}(1+r)\left(\log\frac{e}{(1-r)^2}\right)^{-1}
		\int_{0}^{1}\frac{t}{1+(t-1)r}\log\frac{1+t+(t-1)r}{(1-t)(1-r)}dt\right)\|f\|_{*}.
	\end{align*}
	The above argument establishes the upper bound.

	Now let $f(z)=1$, so that $\|f\|_{\mathcal{B}}=1$. Then
	\begin{align*}
		\|\mathcal{H}\|_{\mathcal{B}\rightarrow   \mathcal{B}_{\log}}&\geq|\mathcal{H}f(0)|+\sup_{z\in\mathbb{D}}(1-|z|^2)\left(\log\frac{e}{(1-|z|)^2}\right)^{-1}|(\mathcal{H}f)'(z)|\\
		&=\left|\int_{0}^{1}f(t)dt\right|+\sup_{z\in\mathbb{D}}(1-|z|^2)\left(\log\frac{e}{(1-|z|)^2}\right)^{-1}\left|\int_{0}^{1}\frac{\phi_t(z)}{1-z}f(\phi_t(z))dt\right|\\
		&\geq1+\sup_{0\leq r<1}(1+r)\left(\log\frac{e}{(1-r)^2}\right)^{-1}\int_{0}^{1}\frac{t}{1+(t-1)r}dt.
	\end{align*}
Next, let $g(z)=\frac{1}{2}\log\frac{1+z}{1-z}$, then $\|g\|_{\mathcal{B}}=1$. We  get that
	\begin{align*}
		&\|\mathcal{H}\|_{\mathcal{B}\rightarrow   \mathcal{B}_{\log}}\geq|\mathcal{H}g(0)|+\sup_{z\in\mathbb{D}}(1-|z|^2)\left(\log\frac{e}{(1-|z|)^2}\right)^{-1}|(\mathcal{H}g)'(z)|\\
		&=\left|\int_{0}^{1}g(t)dt\right|+\sup_{z\in\mathbb{D}}(1-|z|^2)\left(\log\frac{e}{(1-|z|)^2}\right)^{-1}\left|\int_{0}^{1}\frac{\phi_t(z)}{1-z}g(\phi_t(z))dt\right|\\
		&\geq\log2+\frac{1}{2}\sup_{0\leq r<1}(1+r)\left(\log\frac{e}{(1-r)^2}\right)^{-1}
		\int_{0}^{1}\frac{t}{1+(t-1)r}\log\frac{1+t+(t-1)r}{(1-t)(1-r)}dt.
	\end{align*}
	Thus, we establish the desired result.

\end{proof}

	  \begin{lemma}\label{Le3.2}
Let $$A=1+\sup_{0\leq r<1}(1+r)\left(\log\frac{e}{(1-r)^2}\right)^{-1}\int_{0}^{1}\frac{t}{1-(1-t)r}dt,$$
$$B=\log2+\frac{1}{2}\sup_{0\leq r<1}(1+r)\left(\log\frac{e}{(1-r)^2}\right)^{-1} \int_{0}^{1}\frac{t}{1-(1-t)r}\log\frac{1-r+(1+r)t}{(1-r)(1-t)}dt.$$
Then $\max\{A, B\}=A=\displaystyle\frac32$.

\end{lemma}
\begin{proof}

Let
$$
F(r) = (1+r)  \left(\log \frac{e}{(1-r)^2}\right)^{-1} \int_0^1 \frac{t}{1 - (1 - t)r}  dt.
$$
For \( A \), we first have
$$
A = 1 + \sup_{0 \le r < 1} F(r) \ge 1 + F(0) = 1 + \frac{1}{2} = \frac{3}{2}.
$$

Next, we note that
$$
\int_0^1 \frac{t}{1 - (1 - t)r}  dt = \frac{1}{r} + \frac{1 - r}{r^2} \log(1 - r),
$$
and define
$$
f(r) = \frac{1}{r} + \frac{1 - r}{r^2} \log(1 - r) - \frac{1 + r}{2}
= \frac{1 - r}{r^2} \left( \log(1 - r) + \frac{r(r + 2)}{2} \right).
$$
A straightforward calculation shows that \( f(r) \le 0 \). Therefore,
$$
A \le 1 + \sup_{0 \le r < 1}  (1 + r) \left[ \log \frac{e}{(1 - r)^2} \right]^{-1} \cdot \frac{1 + r}{2}
= 1 + \sup_{0 \le r < 1} \frac{(1 + r)^2}{2 \left[ 1 - 2 \log(1 - r) \right]}.
$$
Let
$$
g(r) = \frac{(1 + r)^2}{2 \left[ 1 - 2 \log(1 - r) \right]}.
$$
Then
$$
g'(r) = \frac{-2(1 + r) \left[ r + (1 - r) \log(1 - r) \right]}{(1 - r) \left[ 1 - 2 \log(1 - r) \right]^2}.
$$
We estimate that \( r + (1 - r) \log(1 - r) \ge 0 \), so \( g'(r) \le 0 \), and hence \( g(r) \) is decreasing. It follows that
$ A \le 1 + g(0) = \frac{3}{2}$. We conclude that \( A = \frac{3}{2} \).

To bound $B$, define
	$$h(r)=\displaystyle\int_{0}^{1}\frac{t}{1-(1-t)r}\log\frac{1-r+(1+r)t}{(1-r)(1-t)}dt.$$
We make the change of variable $x=\frac{1-r+(1+r)t}{(1-r)(1-t)}$, so that $t=\frac{(1-r)(x-1)}{1+r+(1-r)x}$.  After computation, we obtain that
\begin{align*}
h(r)&=\int_{1}^{\infty}\frac{x-1}{x+1}\frac{2(1-r)}{\left(1+r+(1-r)x\right)^2}\log xdx\\
&\leq \int_{1}^{\infty}\frac{2(1-r)\log x}{\left(1+r+(1-r)x\right)^2}dx\\
&=-2\int_{1}^{\infty}\log x d(1+r+(1-r)x)^{-1}\\
&=2\int_{1}^{\infty}\frac{1}{1+r+(1-r)x}\frac{1}{x}dx\\
&=\frac{2}{1+r}\log\frac{2}{1-r}.
\end{align*}
Therefore,
\begin{align*}
B
&\leq \log2+\frac{1}{2}\sup_{0\leq r <1}(1+r)\left(\log\frac{e}{(1-r)^2}\right)^{-1}\frac{2}{1+r}\log\frac{2}{1-r}\\
&=\log2+\sup_{0\leq r <1}\frac{\log2-\log(1-r)}{1-2\log(1-r)}\\
&=\log2+\sup_{0\leq r<1}\frac{1}{2}\left(1+\frac{2\log2-1}{1-2\log(1-r)}\right)\\
&=2\log2\\
&<A.
\end{align*}

\end{proof}
Using Lemmas \ref{Le3.1} and \ref{Le3.2}, we can establish the main theorem of this section.
\begin{theorem}
	The norm of the Hilbert matrix operator acting from $\mathcal{B}$ into $\mathcal{B}_{\log}$ satisfies
	\begin{align*}
	\| \mathcal{H}\|_{\mathcal{B}\rightarrow\mathcal{B}_{\log}}= \frac{3}{2}.
	\end{align*}
\end{theorem}

\section{The norm of the Hilbert matrix $\|\mathcal{H}\|_{H^1 \rightarrow H^1_{\log}}$}\label{s3}

In \cite[p.193]{Dia2}, we know that $\mathcal{H}$ is not bounded on $H^{1}$. However, we get that $\mathcal{H}$ is bounded from $H^{1}$ to $H^{1}_{\log}$ in this section. And we estimate the norm of the Hilbert matrix $\|\mathcal{H}\|_{H^1 \rightarrow H^1_{\log}}$.

\begin{lemma}\label{Le2.1}(Hardy's inequality)\cite[p.48]{Dur1}
	If $f(z)=\sum a_nz^n\in H^1$, then$$\sum_{n=0}^{\infty}\frac{|a_n|}{n+1}\leq \pi\|f\|_{H^1}.$$
\end{lemma}

\begin{lemma}\label{Le2.2}\cite{Liu}
	For $z\in\mathbb{D}$ and $c\in\mathbb{R}$, define$$I_c(z):=\frac{1}{2\pi}\int_{0}^{2\pi}\frac{1}{|1-ze^{-i\theta}|^{1+c}}d\theta.$$
	Then the following statements hold.
	
	(1)If $c<0$, then$$1\leq I_c(z)\leq\frac{\Gamma(-c)}{\Gamma^2(\frac{1-c}{2})}.$$
	
	(2)If $c>0$, then$$1\leq(1-|z|^2)^cI_c(z)\leq\frac{\Gamma(c)}{\Gamma^2(\frac{1+c}{2})}.$$
	
	(3)If $c=0$, then$$\frac{1}{\pi}\leq |z|^2(\log\frac{1}{1-|z|^2})^{-1}I_0(z)\leq1.$$
\end{lemma}

\begin{theorem}
	Upper bound estimate for the norm of the Hilbert matrix operator acting from $H^1$ into $H^1_{\log}$ satisfies $$\|\mathcal{H}\|_{H^1 \rightarrow H^1_{\log}}\leq2\pi\log2.$$
\end{theorem}

\begin{proof}
It follows from Lemma 2.1 in \cite{Lan} that for any $f\in H^1$ with $f(z) = \sum_{k=0}^\infty a_kz^k$, we have
	$$\mathcal{H}f(z)=\frac{1}{1-z}F_f(z),$$
	where
	$$F_f(z)=(1-z)\sum_{n=0}^{\infty}\sum_{k=0}^{\infty}\frac{a_k}{n+k+1}z^n.$$
	For $z\in\mathbb{D}$, write
	\begin{align*}
		F_f(z)&=(1-z)\int_0^1\frac{f(t)}{1-tz}dt
		=\sum_{k=0}^{\infty}a_kK_k(z),
	\end{align*}
	where
	\[
	K_k(z)=\int_0^1t^k\frac{1-z}{1-tz}dt.
	\]
	Since, for $0\leq t\leq1$,
	\[
	\sup_{|z|<1}\left|\frac{1-z}{1-tz}\right|=\frac{2}{1+t},
	\]
	we obtain
	\begin{align*}
		|K_k(z)|&\leq 2\int_0^1\frac{t^k}{1+t}dt
		\leq 2\left(\int_0^1t^kdt\right)\left(\int_0^1\frac{dt}{1+t}\right)
		=\frac{2\log2}{k+1},
	\end{align*}
	where the second inequality follows from Chebyshev's integral inequality, since $t^k$ is increasing and $(1+t)^{-1}$ is decreasing on $[0,1]$.
	Consequently,
	$$|\mathcal{H}f(z)|=\left|\frac{1}{1-z}F_f(z)\right|\leq\left|\frac{2\log2}{1-z}\right|\sum_{k=0}^{\infty}\frac{|a_k|}{k+1}.$$
	 By Lemma \ref{Le2.1}, we obtain
 $$|\mathcal{H}f(z)|\leq\left|\frac{2\pi\log2}{1-z}\right|\|f\|_{H^1}.$$
	Applying Lemma \ref{Le2.2}, we  get that
	\begin{align*}
		\|\mathcal{H}f\|_{H^1_{\log}}&=\sup_{0\leq r<1}\left(\log\frac{e}{(1-r)^2}\right)^{-1}M_1(r,\mathcal{H}f)\notag\\
        &\leq\sup_{0\leq r<1}\left(\log\frac{e}{1-r^2}\right)^{-1} M_1(r,\mathcal{H}f)\notag\\
		&\leq2\pi\log2\|f\|_{H^1}\sup_{0\leq r<1}\displaystyle\frac{\displaystyle\frac{1}{2\pi}\int_{0}^{2\pi}\displaystyle\left|\frac{1}{1-re^{i\theta}}\right|d\theta}{\log\frac{e}{1-r^2}}\notag\\
		&\leq2\pi\log2\|f\|_{H^1}\sup_{0\leq r<1}\frac{\log\frac{1}{1-r^2}}{r^2\log\frac{e}{1-r^2}}.
	\end{align*}
	Now, let $x=\log\frac{1}{1-r^2}$. Then
	$$\|\mathcal{H}\|_{H^1 \rightarrow H^1_{\log}}\leq2\pi\log2\sup_{0\leq x<\infty}\frac{xe^x}{(e^x-1)(1+x)}.$$
	Define $g(x)=\frac{xe^x}{(e^x-1)(1+x)}$. Then $$g'(x)=\frac{e^x(e^x-x^2-x-1)}{(e^x-1)^2(1+x)^2}.$$
	Setting $g'(x_0)=0$  implies  $e^{x_0}=1+x_0+x_0^2$.  It follows that
$$g(x_0)=\frac{x_0e^{x_0}}{(e^{x_0}-1)(1+x_0)}=\frac{1+x_0+x_0^2}{(1+x_0)^2}\leq1.$$
Moreover, since
$$\lim\limits_{x\rightarrow0}g(x)=1~~~ \mbox{and}~~~ \lim\limits_{x\rightarrow\infty}g(x)=1,$$
we conclude that $$\|\mathcal{H}\|_{H^1 \rightarrow H^1_{\log}}\leq2\pi\log2.$$
 We  finish the proof of the theorem.
\end{proof}
\begin{theorem}\label{Th2.2}
	Lower bound estimate for the norm of the Hilbert matrix operator acting from $H^1$ into $H^1_{\log}$ satisfies
	\begin{align*}
		\|\mathcal{H}\|_{H^1 \rightarrow H^1_{\log}}\geq\pi.
	\end{align*}
\end{theorem}

\begin{proof}
	 Define
	$$f_\alpha(z)=\frac{1}{(1-z)^\alpha},~~~  0<\alpha<1, ~~~ z\in\mathbb{D}.$$
	It is straightforward to show that $f_\alpha\in H^1$.

The weighted composition operator $T_t$ acting on a function $f_\alpha$ can be expressed as
	$$T_t(f_\alpha)(z)=\frac{((t-1)z+1)^{\alpha-1}}{(1-t)^\alpha}f_\alpha(z),$$
	and it follows that
	$$\mathcal{H}(f_\alpha)(z)=\left(\int_{0}^{1}\frac{((t-1)z+1)^{\alpha-1}}{(1-t)^\alpha}dt\right)f_\alpha(z).$$
	Therefore, we obtain that
	\begin{align*}
       & \|\mathcal{H}\|_{H^1 \rightarrow H^1_{\log}}
        \geq\frac{\|\mathcal{H}(f_\alpha)\|_{H^1_{\log}}}{\|f_\alpha\|_{H^1}}\\
        &=\frac{\sup_{0\leq r<1}\left(\log\frac{e}{(1-r)^2}\right)^{-1}\displaystyle\frac{1}{2\pi}\int_{0}^{2\pi}\left|\int_{0}^{1}\frac{((t-1)re^{i\theta}+1)^{\alpha-1}}{(1-t)^\alpha}dt
        \frac{1}{(1-re^{i\theta})^\alpha}\right|d\theta}{\sup_{0\leq r<1}\displaystyle\frac{1}{2\pi}\int_{0}^{2\pi}\left|\frac{1}{(1-re^{i\theta})^\alpha}\right|d\theta}\\
	\end{align*}
    In terms of the numerator, let $r = 0$ and we get that
    \begin{align*}
    \|\mathcal{H}\|_{H^1 \rightarrow H^1_{\log}}
        &\geq\frac{ \displaystyle\int_{0}^{1}(1-t)^{-\alpha}dt}{\sup_{0\leq r<1}\displaystyle\frac{1}{2\pi}\int_{0}^{2\pi}\left|\frac{1}{(1-re^{i\theta})^\alpha}\right|d\theta}\\
        &=\frac{\frac{1}{1-\alpha}}{\sup_{0\leq r<1}\displaystyle\frac{1}{2\pi}\int_{0}^{2\pi}\left|\frac{1}{(1-re^{i\theta})^\alpha}\right|d\theta}.
    \end{align*}
By Lemma \ref{Le2.2}, we have that
    $$\sup_{0\leq r<1}\frac{1}{2\pi}\int_{0}^{2\pi}\left|\frac{1}{(1-re^{i\theta})^\alpha}\right|d\theta
    \leq\frac{\Gamma(1-\alpha)}{\Gamma^2(\frac{2-\alpha}{2})}.$$
    Consequently,
    $$
        \|\mathcal{H}\|_{H^1 \rightarrow H^1_{\log}}\geq\frac{\Gamma^2(\frac{2-\alpha}{2})}{(1-\alpha)\Gamma(1-\alpha)}    =\frac{\Gamma^2(\frac{2-\alpha}{2})}{\Gamma(2-\alpha)}.
    $$
	Let $\alpha\rightarrow1^-$, we obtain that
    $$
        \lim_{\alpha\rightarrow1^-}\frac{\Gamma^2(\frac{2-\alpha}{2})}{\Gamma(2-\alpha)}       =\lim_{\alpha\rightarrow1^-}\frac{\Gamma(\frac{2-\alpha}{2})\Gamma(\frac{2-\alpha}{2})}{\Gamma(2-\alpha)}=\frac{\pi}{\sin\frac{\pi}{2}}=\pi.
   $$
    We complete the proof.
\end{proof}

\begin{corollary}
	The norm of the Hilbert matrix operator acting from $H^1$ into $H^1_{\log}$ satisfies
	\[
	\pi\leq \|\mathcal{H}\|_{H^1 \rightarrow H^1_{\log}}\leq2\pi\log2.
	\]
\end{corollary}

\begin{conjecture}
\[
	\|\mathcal{H}\|_{H^1 \rightarrow H^1_{\log}}=\pi.
	\]
\end{conjecture}
\section{The norm of the Hilbert matrix $\|\mathcal{H}\|_{H^{\infty} \rightarrow H^{\infty}_{\log}}$}\label{s4}

    While \cite[p.193]{Dia2} shows that $\mathcal{H}$ is not bounded on $H^{\infty}$, we demonstrate in this section that it is in fact bounded as an operator from $H^{\infty}$ to $H^{\infty}_{\log}$. We also compute the exact norm $\|\mathcal{H}\|_{H^{\infty} \rightarrow H^{\infty}_{\log}}$. Before proving this theorem, we also need the following lemma.

\begin{lemma}\label{Le5.1}
Let
\[
  g(x)=\frac{x e^x}{(e^x-1)(2x+1)},\qquad x>0.
\]
Then there exists a unique number \(\gamma\) such that \(g\) is
strictly decreasing on \((0,\gamma)\) and strictly increasing on
\((\gamma,+\infty)\), where
\[
  e^\gamma=1+\gamma+2\gamma^2.
\]
Moreover,
\[
  \sup_{x>0} g(x)=1.
 \]
\end{lemma}

\begin{proof}
Differentiating \(g\) gives
\[
  g'(x)=
  \frac{e^x\bigl(e^x-1-x-2x^2\bigr)}
  {(e^x-1)^2(2x+1)^2}.
\]
For \(x>0\), the denominator is positive and \(e^x>0\), so
\(\operatorname{sgn} g'(x)=\operatorname{sgn} h(x)\), where
\[
  h(x)=e^x-1-x-2x^2.
\]

It remains to show that \(h\) has a unique positive zero.
A direct computation yields
\[
  h'(x)=e^x-1-4x,\qquad h''(x)=e^x-4.
\]
Thus \(h'\) is strictly decreasing on \((0,\ln4)\) and strictly
increasing on \((\ln4,+\infty)\). Since
\(h'(0)=0\) and \(h'(\ln4)=3-4\ln4<0\), we have
\(h'(x)<0\) on \((0,\ln4]\). Also \(h'(x)\to+\infty\) as
\(x\to+\infty\), so there is a unique \(\beta>\ln4\) such that
\(h'(\beta)=0\). Hence \(h\) decreases on \((0,\beta)\) and increases
on \((\beta,+\infty)\). Because \(h(0)=0\) and \(h(x)\to+\infty\)
as \(x\to+\infty\), there exists a unique \(\gamma>\beta\) with
\(h(\gamma)=0\), i.e.
\(e^{\gamma}=1+\gamma+2\gamma^2\). Therefore
\[
  h(x)<0 \quad (0<x<\gamma),\qquad h(x)>0 \quad (x>\gamma),
\]
and consequently
\[
  g'(x)<0 \quad (0<x<\gamma),\qquad g'(x)>0 \quad (x>\gamma).
\]
Thus \(g\) is strictly decreasing on \((0,\gamma)\) and strictly
increasing on \((\gamma,+\infty)\).

Observing that
\[
  \lim_{x\to 0^+} g(x)=1,\qquad
  \lim_{x\to+\infty} g(x)=\frac12,
\]
and using the monotonicity of g, we have
$g(x)<1$ for all $x>0$. Hence
\[
  \sup_{x>0} g(x)=1.
\]
\end{proof}

\begin{theorem}
	The norm of the Hilbert matrix operator acting from $H^\infty$ into $H^\infty_{\log}$ satisfies$$\|\mathcal{H}\|_{H^\infty \rightarrow H^\infty_{\log}}=1.$$
\end{theorem}
\begin{proof}
	Let $f\in H^\infty$. Then we have
	\begin{align*}
		\|\mathcal{H}(f)\|_{H^\infty_{\log}}&=\sup_{z\in\mathbb{D}}\left(\log\frac{e}{(1-|z|)^2}\right)^{-1}\left|\int_{0}^{1}\frac{f(t)}{1-tz}dt\right|\\
		&\leq\sup_{z\in\mathbb{D}}\left(\log\frac{e}{(1-|z|)^2}\right)^{-1}\int_{0}^{1}\left|\frac{1}{1-tz}\right|dt\|f\|_{H^\infty}\\
		&\leq\sup_{0\leq r<1}\left(\log\frac{e}{(1-r)^2}\right)^{-1}\int_{0}^{1}\frac{1}{1-tr}dt\|f\|_{H^\infty}\\
		&=\sup_{0< r<1}\frac{\displaystyle\frac{1}{r}\log\frac{1}{1-r}}{\log\frac{e}{(1-r)^2}}\|f\|_{H^\infty}.
	\end{align*}
For the lower bound,  take $f=1\in H^\infty$. Then
	$$\|f\|_{H^\infty}=1, \mathcal{H}(1)(z)=\frac{1}{z}\log\frac{1}{1-z}.$$
	Thus,
	\begin{align*}
		\|\mathcal{H}\|_{H^\infty \rightarrow H^\infty_{\log}}
        &\geq\sup_{z\in\mathbb{D}}\left(\log\frac{e}{(1-|z|)^2}\right)^{-1}
        \left|\frac{1}{z}\log\frac{1}{1-z}\right|\\
		&\geq\sup_{0< r<1}\frac{\displaystyle\frac{1}{r}\log\frac{1}{1-r}}{\log\frac{e}{(1-r)^2}}.
	\end{align*}
    We conclude that
    $$\|\mathcal{H}\|_{H^\infty \rightarrow H^\infty_{\log}}=\sup_{0< r<1}\frac{\displaystyle\frac{1}{r}\log\displaystyle\frac{1}{1-r}}{\log\frac{e}{(1-r)^2}}.$$
  Then, making the change of variables  $x=\log\frac{1}{1-r}$, we obtain
    $$\|\mathcal{H}\|_{H^\infty \rightarrow H^\infty_{\log}}=\sup_{x>0}\frac{xe^{x}}{(e^{x}-1)(2x+1)}.$$
    By  Lemma \ref{Le5.1}, we have
     $$\|\mathcal{H}\|_{H^\infty \rightarrow H^\infty_{\log}}=1.$$
This completes the proof.

\end{proof}
\section{The norm of the Hilbert matrix $\|\mathcal{H}\|_{\mathcal{B}^\alpha \rightarrow \mathcal{B}^\alpha_{\log}}$  with $1<\alpha<2$}\label{s6}
\begin{theorem}
	For $1<\alpha<2$, we have that
	$$\|\mathcal{H}\|_{\mathcal{B}^\alpha\rightarrow \mathcal{B}^\alpha_{\log}}
    \geq\frac{1}{2(\alpha-1)}\int_{0}^{1}(1-t^2)^{1-\alpha}dt +\frac{3\alpha-5}{4(\alpha-1)(2-\alpha)},$$
	and $\mathcal{H}$ cannot  map $\mathcal{B}^\alpha$ to $\mathcal{B}^\alpha_{\log}$ for $0<\alpha<1$ and $\alpha\geq2$.
\end{theorem}
\begin{proof}
	Let $\alpha\not=1$ and $z\in\mathbb{D}$. Define$$f_\alpha(z)=\frac{1}{2(\alpha-1)}(1-z^2)^{1-\alpha}-\frac{1}{2(\alpha-1)}$$
	On one hand, we have the estimate
	$$\|f_\alpha\|_{\mathcal{B}^\alpha}=\sup_{z\in\mathbb{D}}\frac{|z|(1-|z|^2)^\alpha}{|1-z^2|^\alpha}\leq\sup_{z\in\mathbb{D}}|z|=1.$$
	On the other hand, for $r\in(0,1)$ it holds that $\lim_{r\rightarrow1^-}|f'_\alpha(r)|(1-r^2)^\alpha=1$, and we obtain $\|f_\alpha\|_{\mathcal{B}^\alpha}=1.$
	
	According to (\ref{eq3.2}) , we obtain that
	\begin{align*}
		&\|\mathcal{H}\|_{\mathcal{B}^\alpha\rightarrow\mathcal{B}^\alpha_{\log}}\geq\frac{\|\mathcal{H}f_\alpha\|_{\mathcal{B}^\alpha_{\log}}}
{\|f_\alpha\|_{\mathcal{B}^\alpha}}=|\mathcal{H}f_\alpha(0)|+\sup_{z\in\mathbb{D}}\left(\log\frac{e}{(1-|z|)^2}\right)^{-1}(1-|z|^2)^\alpha|(\mathcal{H}f_\alpha)'(z)|\\
		&=\left|\frac{1}{2(\alpha-1)}\left(\int_{0}^{1}(1-t^2)^{1-\alpha}dt-1\right)\right|\\
		&+\frac{1}{2|\alpha-1|}\sup_{z\in\mathbb{D}}\left(\log\frac{e}{(1-|z|)^2}\right)^{-1}(1-|z|^2)^\alpha\left|\int_{0}^{1}\frac{t(1-\phi_t(z)^2)^{1-\alpha}-t}{\left[(t-1)z+1\right](1-z)}dt\right|\\
		&\geq\left|\frac{1}{2(\alpha-1)}\int_{0}^{1}(1-t^2)^{1-\alpha}dt-\frac{1}{2(\alpha-1)}\right|+\frac{1}{2|\alpha-1|}\sup_{z\in\mathbb{D}}\left(\log\frac{e}{(1-|z|)^2}\right)^{-1}(1-|z|^2)^\alpha\\
		&\times\left|\int_{0}^{1}\frac{t(1-t)^{1-\alpha}\left[(t-1)z+1+t\right]^{1-\alpha}}{(1-z)^\alpha\left[(t-1)z+1\right]^{3-2\alpha}}-\frac{t}{\left[(t-1)z+1\right](1-z)}dt\right|.
	\end{align*}
	For $0<\alpha<1$, we have that
	\begin{align*}
		&\frac{1}{2(1-\alpha)}\sup_{z\in\mathbb{D}}\left(\log\frac{e}{(1-|z|)^2}\right)^{-1}(1-|z|^2)^\alpha\\
&\times \left|\int_{0}^{1}\frac{t(1-t)^{1-\alpha}\left[(t-1)z+1+t\right]^{1-\alpha}}{(1-z)^\alpha\left[(t-1)z+1\right]^{3-2\alpha}}-\frac{t}{\left[(t-1)z+1\right](1-z)}dt\right|\\
		&\geq\frac{1}{2(1-\alpha)}\lim\limits_{r\rightarrow1^-}\left(\log\frac{e}{(1-r)^2}\right)^{-1}(1-r^2)^\alpha\\
&\times\left|\int_{0}^{1}\frac{t}{\left[(t-1)r+1\right](1-r)}-\frac{t(1-t)^{1-\alpha}\left[(t-1)r+1+t\right]^{1-\alpha}}{(1-r)^\alpha\left[(t-1)r+1\right]^{3-2\alpha}}dt\right|\\
		&=\infty.
	\end{align*}
Hence, $\mathcal{H}$ cannot map $\mathcal{B}^\alpha$ into $\mathcal{B}^\alpha_{\log}$ for $0 < \alpha < 1$.	

	For $\alpha\geq2$, the integral $\left|\displaystyle\int_{0}^{1}\frac{1}{2(\alpha-1)}(1-t^2)^{1-\alpha}dt-\frac{1}{2(\alpha-1)}\right|$ diverges,
 so  $\mathcal{H}$ cannot  map $\mathcal{B}^\alpha$ to $\mathcal{B}^\alpha_{\log}$ for $\alpha\geq2$.
	
	 For $1<\alpha<2$, we obtain that
	 \begin{align*}
	 	&\|\mathcal{H}\|_{\mathcal{B}^\alpha\rightarrow\mathcal{B}^\alpha_{\log}}\\
	 	&\geq\frac{1}{2(\alpha-1)}\int_{0}^{1}(1-t^2)^{1-\alpha}dt-\frac{1}{2(\alpha-1)}+\frac{1}{2(\alpha-1)}\sup_{0\leq r<1}\left(\log\frac{e}{(1-r)^2}\right)^{-1}(1-r^2)^\alpha\\ &\times\left|\int_{0}^{1}\frac{t(1-t)^{1-\alpha}\left[(t-1)r+1+t\right]^{1-\alpha}}{(1-r)^\alpha\left[(t-1)r+1\right]^{3-2\alpha}}-\frac{t}{\left[(t-1)r+1\right](1-r)}dt\right|\\
	 	&\geq\frac{1}{2(\alpha-1)}\int_{0}^{1}(1-t^2)^{1-\alpha}dt-\frac{1}{2(\alpha-1)}+\frac{1}{2(\alpha-1)}\lim_{r\rightarrow0}\left(\log\frac{e}{(1-r)^2}\right)^{-1}(1-r^2)^\alpha\\
	 	&\times\left|\int_{0}^{1}\frac{t(1-t)^{1-\alpha}\left[(t-1)r+1+t\right]^{1-\alpha}}{(1-r)^\alpha\left[(t-1)r+1\right]^{3-2\alpha}}-\frac{t}{\left[(t-1)r+1\right](1-r)}dt\right|\\
	 	&=\frac{1}{2(\alpha-1)}\int_{0}^{1}(1-t^2)^{1-\alpha}dt+\frac{3\alpha-5}{4(\alpha-1)(2-\alpha)}.
	\end{align*}
	This completes the proof.
\end{proof}
Next, we consider the upper bound of the norm estimate for Hilbert matrix operator acting from $\mathcal{B}^\alpha$ into $\mathcal{B}_{\log}^\alpha$ for $1<\alpha<2$.
\begin{theorem}
	For $1<\alpha<2$, the norm of the Hilbert matrix operator acting  from $\mathcal{B}^\alpha$ into $\mathcal{B}_{\log}^\alpha$  satisfies $$\|\mathcal{H}\|_{\mathcal{B}^\alpha\rightarrow\mathcal{B}^\alpha_{\log}}\leq\frac{\pi}{\sin(\alpha-1)\pi}+\frac{1}{2-\alpha}.$$
\end{theorem}
\begin{proof}  Suppose $f\in\mathcal{B}^\alpha$ for $1<\alpha<2$ and $z\in\mathbb{D}$. 	From (\ref{eq3.2}) and (\ref{eq3.4}) we can deduce that
	\begin{align*}
		&\|\mathcal{H}f\|_{\mathcal{B}^\alpha_{\log}}=|\mathcal{H}f(0)|+\sup_{z\in\mathbb{D}}\left(\log\frac{e}{(1-|z|)^2}\right)^{-1}(1-|z|^2)^\alpha|(\mathcal{H}f)'(z)|\\
		&=\left|\int_{0}^{1}f(t)dt\right|+\sup_{z\in\mathbb{D}}\left(\log\frac{e}{(1-|z|)^2}\right)^{-1}(1-|z|^2)^\alpha\left|\int_{0}^{1}\frac{\phi_t(z)}{1-z}f(\phi_t(z))dt\right|\\
		&\leq\int_{0}^{1}|f(t)|dt+\sup_{z\in\mathbb{D}}\left(\log\frac{e}{(1-|z|)^2}\right)^{-1}(1-|z|^2)^\alpha\int_{0}^{1}\frac{|\phi_t(z)|}{1-|z|}|f(\phi_t(z))|dt\\
		&\leq\int_{0}^{1}\left(\frac{(1-t)^{1-\alpha}-1}{\alpha-1}\|f\|_{\alpha*}+|f(0)|\right)dt\\
		&+\sup_{z\in\mathbb{D}}\left(\log\frac{e}{(1-|z|)^2}\right)^{-1}(1+|z|)^\alpha(1-|z|)^{\alpha-1}\\
&\times\int_{0}^{1}|\phi_t(z)|\bigg[\frac{(1-\phi_t(z))^{1-\alpha}-1}{\alpha-1}\|f\|_{\alpha*}+|f(0)|\bigg]dt.
	\end{align*}

	Since $|\phi_t(z)|=|\frac{t}{1+(t-1)z}|\leq\frac{t}{1+(t-1)|z|}$ and $1<\alpha<2$, we find that
	
	\begin{align*}
		\|\mathcal{H}f\|_{\mathcal{B}^\alpha_{\log}}
&\leq\frac{1}{2-\alpha}\|f\|_{\alpha*}+|f(0)|\\ &+\sup_{z\in\mathbb{D}}\left(\log\frac{e}{(1-|z|)^2}\right)^{-1}(1+|z|)^\alpha(1-|z|)^{\alpha-1}\\
&\times\int_{0}^{1}\frac{t}{1+(t-1)|z|}\left(\frac{(1-t)^{1-\alpha}(1-|z|)^{1-\alpha}}{(\alpha-1)(1+(t-1)|z|)^{1-\alpha}}-\frac{1}{\alpha-1}\right)\|f\|_{\alpha*}dt\\ &+\sup_{z\in\mathbb{D}}\left(\log\frac{e}{(1-|z|)^2}\right)^{-1}(1+|z|)^\alpha(1-|z|)^{\alpha-1}\int_{0}^{1}\frac{t}{1+(t-1)|z|}|f(0)|dt\\
&=\frac{1}{2-\alpha}\|f\|_{\alpha*}+|f(0)|+\sup_{0\leq r<1}\left(\log\frac{e}{(1-r)^2}\right)^{-1}(1+r)^\alpha\\
&\times\int_{0}^{1}\frac{t}{1+(t-1)r}\left(\frac{(1-t)^{1-\alpha}}{(\alpha-1)(1+(t-1)r)^{1-\alpha}}-\frac{(1-r)^{\alpha-1}}{\alpha-1}\right)\|f\|_{\alpha*}dt\\
&+\sup_{0\leq r<1}\left(\log\frac{e}{(1-r)^2}\right)^{-1}(1+r)^\alpha(1-r)^{\alpha-1}
\int_{0}^{1}\frac{t}{1+(t-1)r}|f(0)|dt\\
&=\bigg[\frac{1}{2-\alpha}+\sup_{0\leq r<1}\left(\log \frac{e}{(1-r)^2}\right)^{-1}(1+r)^\alpha\\
&\times\int_{0}^{1}\frac{t}{1+(t-1)r}\left(\frac{(1-t)^{1-\alpha}}{(\alpha-1)(1+(t-1)r)^{1-\alpha}}-\frac{(1-r)^{\alpha-1}}{\alpha-1}\right)dt\bigg]\|f\|_{\alpha*}\\
&+\left(1+\sup_{0\leq r<1}\left(\log\frac{e}{(1-r)^2}\right)^{-1}(1+r)^\alpha(1-r)^{\alpha-1}\int_{0}^{1}\frac{t}{1+(t-1)r}dt\right)|f(0)|.
	\end{align*}
	Since the function  $\displaystyle\frac{1}{(1+(t-1)r)^{2-\alpha}}$ is increasing with respect to
 $r$, we obtain that
\begin{align}\label{eq3.7}
		&\sup_{0\leq r<1}\left(\log\frac{e}{(1-r)^2}\right)^{-1}\int_{0}^{1}\frac{t(1+r)^\alpha}{1+(t-1)r}\left(\frac{(1-t)^{1-\alpha}}{(\alpha-1)(1+(t-1)r)^{1-\alpha}}-\frac{(1-r)^{\alpha-1}}{\alpha-1}\right)dt  \notag\\
        &\leq\sup_{0\leq r<1}\left(\log\frac{e}{(1-r)^2}\right)^{-1}
\frac{(1+r)^\alpha}{\alpha-1}\int_{0}^{1}\frac{t(1-t)^{1-\alpha}}{(1+(t-1)r)^{2-\alpha}}dt  \notag\\
                &\leq \sup_{0\leq r<1}\left(\log\frac{e}{(1-r)^2}\right)^{-1}
\frac{(1+r)^\alpha}{\alpha-1}\int_{0}^{1}t^{\alpha-1}(1-t)^{1-\alpha}dt   \notag\\
                      &=\frac{\pi}{\sin(\alpha-1)\pi}\sup_{0\leq r<1}\left(\log \frac{e}{(1-r)^2}\right)^{-1}
                      (1+r)^\alpha\leq \frac{\pi}{\sin(\alpha-1)\pi},
    \end{align}
where we use the fact that $\displaystyle\int_{0}^{1}t^{\alpha-1}(1-t)^{1-\alpha}dt=\frac{(\alpha-1)\pi}{\sin(\alpha-1)\pi}$ and
$\displaystyle\sup_{0\leq r<1}\left(\log\frac{e}{(1-r)^2}\right)^{-1}(1+r)^\alpha\leq 1$.
 Then we have  that
\begin{align*} \|\mathcal{H}f\|_{\mathcal{B}^\alpha_{\log}}&\leq\left(\frac{1}{2-\alpha}+\frac{\pi}{\sin(\alpha-1)\pi}\right)\|f\|_{\alpha*}\\
		&+\left(1+\sup_{0\leq r<1}\left(\log\frac{e}{(1-r)^2}\right)^{-1}(1+r)^\alpha(1-r)^{\alpha-1}\int_{0}^{1}\frac{t}{1+(t-1)r}dt \right)|f(0)|.
\end{align*}
On the other hand, for $1<\alpha<2$, we note that
   $$(1-r)^{\alpha-1}\leq \frac{(1-t)^{1-\alpha}}{(\alpha-1)(1-r)^{1-\alpha}}
\leq \frac{(1-t)^{1-\alpha}}{(\alpha-1)(1+(t-1)r)^{1-\alpha}},$$
from which, together with (\ref{eq3.7}), we deduce
    \begin{align*}
    &\sup_{0\leq r<1}\left(\log\frac{e}{(1-r)^2}\right)^{-1}(1+r)^\alpha(1-r)^{\alpha-1}
\int_{0}^{1}\frac{t}{1+(t-1)r}dt\\
    &\leq \sup_{0\leq r<1}\left(\log\frac{e}{(1-r)^2}\right)^{-1}
\frac{(1+r)^\alpha}{\alpha-1}\int_{0}^{1}\frac{t(1-t)^{1-\alpha}}{(1+(t-1)r)^{2-\alpha}}dt\\
    &\leq \frac{\pi}{\sin(\alpha-1)\pi}.
    \end{align*}
Moreover, it is obvious that $1< \frac{1}{2-\alpha}$ for $1<\alpha<2$.
We thus conclude that
$$
1+\sup_{0\leq r<1}\left(\log\frac{e}{(1-r)^2}\right)^{-1}(1+r)^\alpha(1-r)^{\alpha-1}
\int_{0}^{1}\frac{t}{1+(t-1)r}dt
\leq \frac{1}{2-\alpha}+\frac{\pi}{\sin(\alpha-1)\pi}.
$$
	Therefore,
$$\|\mathcal{H}\|_{\mathcal{B}^\alpha\rightarrow\mathcal{B}^\alpha_{\log}}\leq\frac{\pi}{\sin(\alpha-1)\pi}+\frac{1}{2-\alpha}.$$
We complete the proof.
\end{proof}
\begin{corollary}
	The Hilbert matrix operator $\mathcal{H}$ maps $\mathcal{B}^\alpha$ boundedly into $\mathcal{B}^\alpha_{\log}$ if and only if
	\[
	1\leq\alpha<2.
	\]
\end{corollary}
\begin{proof}
	The case $\alpha=1$ follows from Lemmas \ref{Le3.1} and \ref{Le3.2}. The case $1<\alpha<2$ follows from the preceding theorem. The non-boundedness for $0<\alpha<1$ and $\alpha\geq2$ was proved in the first theorem of this section.
\end{proof}

\noindent{\bf Acknowledgments}

\noindent   The research was supported by Zhejiang Provincial Natural Science Foundation (Grant No.
 LY23A010003) and National Natural Science Foundation of China (Grant No. 11671357).

\end{document}